\documentclass[12pt]{svjour3}
\usepackage{amssymb}
\usepackage{amsmath}
\usepackage{color}
\usepackage{comment}
\usepackage{enumitem}
\usepackage{graphicx}

\begin{document}

\def\x{{y}}
\def\sp{{y}}
\def\bW{{\bf W}}
\def\bx{{\bf x}}
\def\by{{\bf y}}
\def\cF{{\mathcal{F}}}
\def\ee{\varepsilon}
\def\T{{T}}

\newcommand{\removableFootnote}[1]{\footnote{#1}}

\newcommand{\fig}[1]{Fig.~\ref{fig:#1}}

\title{Non-Filippov dynamics arising from the smoothing of nonsmooth systems, and its robustness to noise}
\titlerunning{Non-Filippov dynamics}

\author{M.R. Jeffrey \and D.J.W. Simpson}
\institute{M.R. Jeffrey \at 
Department of Engineering Mathematics\\
University of Bristol, Bristol, UK
\and
D.J.W. Simpson \at
Institute of Fundamental Sciences\\
Massey University\\
Palmerston North, New Zealand}

\date{\today}
\maketitle

\begin{abstract}

Switch-like behaviour in dynamical systems may be modelled by highly nonlinear functions, such as Hill functions or sigmoid functions, or alternatively by piecewise-smooth functions, such as step functions. 
Consistent modelling requires that piecewise-smooth and smooth dynamical systems have similar dynamics, but the conditions for such similarity are not well understood. Here we show that by smoothing out a piecewise-smooth system one may obtain dynamics that is inconsistent with the accepted wisdom --- so-called {\it Filippov} dynamics --- at a discontinuity, even in the piecewise-smooth limit. 
By subjecting the system to white noise, we show that these discrepancies can be understood in terms of potential wells that allow solutions to dwell at the discontinuity for long times. Moreover we show that spurious dynamics will revert to Filippov dynamics, with a small degree of stochasticity,
when the noise magnitude is sufficiently large compared to the order of smoothing.
We apply the results to a model of a dry-friction oscillator, where spurious dynamics (inconsistent with Filippov's convention or with Coulomb's model of friction) can account for different coefficients of static and kinetic friction, 
but under sufficient noise the system reverts to dynamics consistent with Filippov's convention (and with Coulomb-like friction).

\end{abstract}

\section{Introduction}
\label{sec:INTRO}
\setcounter{equation}{0}

Systems of piecewise-smooth differential equations 
are used to model diverse phenomena involving switch-like and impact dynamics throughout engineering, biology, and physics;
see for example \cite{bc08,krg03,l04,physDspecial} and references therein.
They consist of ordinary differential equations that are smooth except at certain {\it switching surfaces}, where the presence of a discontinuity permits a variety of intricate dynamical behaviour and bifurcations that are not possible in smooth systems.
The study of how dynamical systems theory can incorporate discontinuities has been the subject of much recent research (see e.g. \cite{bc08,cdhj11special} and references therein).
The extent to which any piecewise-smooth system can be considered as an approximation to a smooth system, or vice versa, remains an open question.
Here we show that smoothing a discontinuity can result in unexpected dynamics, which can be further understood by considering the response of the system to noise. 

The problem of how to model switch-like behaviour in dynamical systems is one of considerable complexity. As well as certain unobvious difficulties in relating smooth and piecewise-smooth approximations, physical applications show that switches are particularly prone to effects of hysteresis, time-delay, and noise. In the seemingly elementary mechanics problem of dry-friction between rigid bodies, there remain numerous viewpoints as to the optimum way to improve upon Coulomb's simple discontinuous contact laws, see for example \cite{Bliman95easy-to-userealistic,popp,krim,wojewoda08}. Akay considered vibrations due to friction and their effect on friction dynamics \cite{akay02,akay99}. 
The effects of surface roughness and stochastic motions on the friction force have been studied experimentally \cite{feldmann12,bot09,persson05}. The problem has also been studied theoretically, with de Gennes proposing the study of dry-friction on surfaces subjected to external white noise vibrations \cite{gennes05}, and this has been extended via a path integral approach to account for basic stick-slip motion \cite{baule-pathintegral,baule-stickslip}.
The likely sources of noise in electronics, and ways to minimise its effects through regularisation, are of particular interest in the control literature \cite{db02,slot91,u92}. 
While many piecewise-smooth approaches to switching dynamics follow Filippov's approach \cite{bc08,f88,krg03}, attempts to relate their properties to smooth models \cite{bps01,t07} raise the question of how general the approach is, and while alternative views have certainly been expressed \cite{aizerman12,hajek1,hermes68,j13error}, a definitive answer to how closely smooth and piecewise-smooth systems approximate each other remains lacking. Here we show how dynamics outside the Filippov convention for piecewise-smooth systems can arise, and how it can be understood in the presence of perturbations due to noise. 

In short, the standard `Filippov' approach involves restricting the form of a set of differential equations {\it at} a discontinuity, to being a linear interpolation of its two forms immediately {\it either side} of the discontinuity. As Filippov showed \cite{f88}, the equations then have solutions that are sufficiently unique to represent a deterministic dynamical system (with the exception of certain singularities \cite{j11prl}). The physical significance of the Filippov convention has been proven repeatedly in application to electronic control, stick-slip and other mechanical behaviours, and various biological models (see e.g. \cite{bc08,slot91,u92,physDspecial} and references therein). In \cite{j13error} it was shown that spurious dynamics, which lies outside the scope of the Filippov approach, can not only be introduced analytically into a deterministic model, but persists when a system is smoothed out or subjected to random perturbations. 

The aim here is to study the effect of random perturbations more closely and more rigorously, by focussing on the effect of white noise. 
We study what happens to dynamics near a discontinuity when the governing equations are smoothed out and when noise is added. These two perturbations are motivated by various practical and theoretical considerations. 
In control applications, the {\it regularisation} of a discontinuous switch by smoothing
is utilised to avoid wear and instability \cite{slot91}.
In mechanics it is often unclear whether a smooth or discontinuous model is more appropriate \cite{j13error,krim}. 
Mathematically, though one can smooth out a discontinuity and show topological equivalence to classes of slow-fast systems \cite{t07}, this does not take into account spurious dynamics beyond the Filippov model. 
On the problem of noise, it has recently been shown that stochastic solutions can reduce to the deterministic Filippov solution in the zero noise limit \cite{boq09,simpson12}, again assuming a priori the absence of effects that defy the Filippov convention. 

The present paper shows that spurious dynamics predicted in \cite{j13error}
can be understood through stochastic dynamics.
In the presence of noise, smoothing out the discontinuity can create local potential wells where the solutions may dwell for extremely long times, even in regions where (\ref{eq:genFilippov}) predicts an immediate escape. 
We show that the robustness of such spurious dynamics depends on the relative length-scales of the perturbations due to noise and due to smoothing.

The systems of interest are defined in \S\ref{sec:SYSTEMS}, and for convenience our main analytical results are then summarised in \S\ref{sec:RESULTS}. In \S\ref{sec:SCALE} we obtain a local approximation, from which one-dimensional Fokker-Planck equations are derived in \S\ref{sec:STOCH}. The implied escape times and probability density functions used to describe sliding and crossing dynamics, respectively, are studied asymptotically in \S\ref{sec:CALC}, with some details of calculations included in an Appendix. In \S\ref{sec:FRICTION} we apply these results to a dry-friction oscillator, showing that spurious sliding dynamics can be used to model disparity of the static and kinetic friction coefficients, quantifying its robustness to noise. Closing remarks, including suggestions for further work, are made in \S\ref{sec:CONC}.

\section{Description of the problem}\label{sec:SYSTEMS}

To investigate the dynamics of a piecewise-smooth system whose switching surface is a manifold in $\mathbb{R}^n$, we first note that with a suitable choice of coordinates $\bx \in \mathbb{R}^n$, the switching surface is a coordinate plane, say $\x=0$ (where $\x$ is a component of the vector $\bx$).
The system may then be written as
\begin{equation}
\dot{\bx} = 
\left\{ \begin{array}{lc}
f^-(\bx) \;, & \x < 0 \\
f^+(\bx) \;, & \x > 0
\end{array} \right. \;,
\label{eq:genFilippov}
\end{equation}
and we assume $f^-$ and $f^+$ are smooth functions. In the remainder of this section we summarise the approach for solving (\ref{eq:genFilippov}), and state our main results concerning the effect of perturbations. 

\begin{figure}[h!]
\begin{center}
\includegraphics[width=0.7\textwidth]{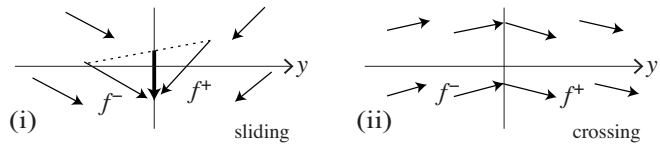}
\caption{\sf\small 
Filippov dynamics in the system (\ref{eq:genFilippov}) near $\x=0$ in two dimensions.
Panel (i) shows an attracting sliding region,
where the vector field points toward the switching surface $\x=0$ from both sides, and on $\x=0$ the flow follows the {\it sliding vector field}, defined as the tangent vector to $\x=0$ given by a convex combination of $f^+$ and $f^-$ (bold arrow). 
Panel (ii) shows a crossing region,
where the vector field points through $\x=0$ and the flow crosses instantaneously.
\label{fig:fil}
}
\end{center}
\end{figure}

As described by Filippov \cite{f88},
there are two important, fundamental scenarios for the dynamics near the switching surface: 
{\it sliding} and {\it crossing}, illustrated in \fig{fil}.
These two behaviours form a foundation for the study of piecewise-smooth dynamical systems.
Other scenarios, such as a repelling form of sliding (see \cite{f88,j11prl,js09,jh09}),
or tangencies between the vector field and the switching surface (see \cite{bc08,f88,krg03,simpson13}),
are not considered here. 

Of the two scenarios in \fig{fil}, sliding is usually assumed to be the most novel.
When both of the vector fields $f^\pm$ point toward the switching surface, they constrain solutions of (\ref{eq:genFilippov}) to evolve along $y=0$.
The velocity of this evolution is conventionally given by a convex combination of $f^+$ and $f^-$.
The diverse applications of sliding include
frictional stick-slip models \cite{bc08,cdhj11special},
sliding-mode control \cite{u77},
and prey-switching in animal predation \cite{dgr07,krg03,piltz13}.
The crossing scenario occurs when the components of $f^-$ and $f^+$
normal to $\x=0$ have the same sign, as in \fig{fil}(ii), then solutions of (\ref{eq:genFilippov}) are assumed to cross $\x=0$ instantaneously. 

\begin{figure}[h!]
\begin{center}
\includegraphics[width=\textwidth]{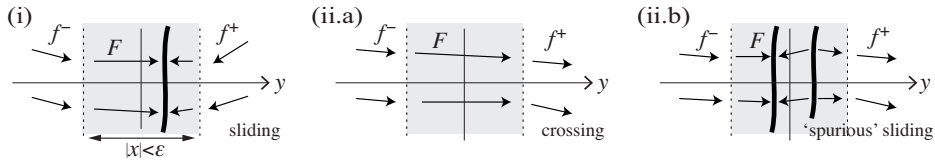}
\caption{\sf\small 
Sketches of the system (\ref{eq:genPerturbed}) in two dimensions.
Panel (i) shows a smoothing \fig{fil}(i), where a bold curve shows the locus of sliding-like dynamics.
Panels (ii.a) and (ii.b) show different smoothings of \fig{fil}(ii). In panel (ii.a) the smoothing function $F$ is linear, resulting in crossing-like dynamics consistent with Filippov dynamics. 
In panel (ii.b), the first component of $F$ changes sign twice within the region $|\x|<\ee$, creating spurious sliding-like dynamics (of both attracting and repelling types). 
\label{fig:flows}
}
\end{center}
\end{figure}

We then consider two perturbations of (\ref{eq:genFilippov})
motivated by general modelling considerations.
First we make (\ref{eq:genFilippov}) continuous by supposing that
the jump between $f^+$ and $f^-$ takes place within the region $|\x|\le\ee$ for some small $\ee>0$.
Specifically we write
\begin{equation}
\dot{\bx} = 
\left\{ \begin{array}{lc}
f^-(\bx) \;, & \x \le -\ee\;, \\
F \left( \frac{\x}{\ee}, \bx \right) \;, & -\ee < \x < \ee\;, \\
f^+(\bx) \;, & \x \ge \ee\;,
\end{array} \right. 
\label{eq:genPerturbed}
\end{equation}
where $F(\frac{\x}{\ee},\bx)$ is a smooth
function of $\frac{\x}{\ee}$ and $\bx$.
(For the asymptotic expansions that we employ below, it is useful to explicitly keep track of the ratio $\frac{\x}{\ee}$,
separately to the occurrence of $\x$ inside the vector $\bx$.)
Continuity of the right-hand side of (\ref{eq:genPerturbed}) requires
\begin{equation}
\label{eq:cont}
\lim_{\x\rightarrow+\ee}F\left(\frac{\x}{\ee},\bx\right)=\lim_{\x\rightarrow+\ee} f^+(\bx)\quad\;\;\&\quad
\lim_{\x\rightarrow-\ee}F\left(\frac{\x}{\ee},\bx\right)=\lim_{\x\rightarrow-\ee} f^-(\bx)\;.
\end{equation}
%
%
\fig{flows} illustrates the smoothed system (\ref{eq:genPerturbed}), and corresponds to the smoothing of \fig{fil} with simple choices of $F$. From the sliding scenario of \fig{fil}(i), the continuity condition (\ref{eq:cont}) demands that the $\x$-component of $F$ changes direction inside the region $|\x|<\ee$, giving \fig{flows}(i).
From the crossing scenario of \fig{fil}(ii), the continuity condition implies that
the $\x$-component of $F$ has the same sign at $\x=-\ee$ and $\x=\ee$. This sign may be constant throughout $|\x|<\ee$, as in \fig{flows}(ii.a), or it may change within $|\x|<\ee$, as in \fig{flows}(ii.b). The two cases lead to qualitatively different dynamics, the former resembling the crossing of the piecewise-smooth model in \fig{fil}(i), while the latter exhibits a form of sliding that is inconsistent with Filippov dynamics. In general, the $\x$ component of $F$ may change sign many times in $|\x|<\ee$, but the simplest cases given in \fig{flows} are sufficient to initiate a study of spurious dynamics and its robustness to noise. These scenarios are also sufficient to investigate a novel cause of {\it stiction} in the dry-friction model of \S\ref{sec:FRICTION}.

The second perturbation we consider is the addition of noise, in the form
\begin{equation}
d\bx(t) = 
\left\{ \begin{array}{lc}
f^-(\bx(t)) \;, & \x(t) \le -\ee \\
F \left( \frac{\x(t)}{\ee}, \bx(t) \right) \;, & -\ee < \x(t) < \ee \\
f^+(\bx(t)) \;, & \x(t) \ge \ee
\end{array} \right\} \,dt + \kappa D(\bx(t)) \,d\bW(t) \;,
\label{eq:genSDE}
\end{equation}
where $\bW(t)$ represents a standard vector Brownian motion,
$0 < \kappa \ll 1$ is the noise amplitude,
and $D$ is an $n \times n$ matrix that represents the strength of the noise in different directions.
This models the effect of random external white noise, assuming it is uniform in $\bx$ and is not, for example, skewed or amplified near $\x=0$.
Alternative formulations, such as coloured noise or an $\bx$-dependent noise amplitude, would be interesting extensions for future work.

As we shall see, if we consider (\ref{eq:genSDE}) with $\kappa$ sufficiently large relative to $\ee$,
the effect of noise tends to push solutions outside the discontinuity neighbourhood $|\x|<\ee$, and hence beyond the influence of $F$.
For both sliding and crossing we then recover Filippov's solution, albeit with a small degree of stochasticity. These results are consistent in outcome, if different in set-up, to \cite{simpson12}.
However, if $\kappa$ is relatively small and we consider a nonlinear smoothing function $F$,
the dynamics of (\ref{eq:genSDE}) may be vastly different to that predicted under Filippov's convention.
For instance for \fig{flows}(ii.b), solutions to (\ref{eq:genPerturbed}) are trapped near $\x\approx0$ for a long period of time, whereas Filippov's convention suggests that solutions cross the region $|\x|<\ee$ in a small time of order $\ee$.
This extra time implies that a crossing region may exhibit unexpected dynamics.
However, if $\kappa$ is sufficiently large relative to $\ee$ then noise can dominate these spurious smoothing effects, and then Filippov dynamics is restored. 

The purpose of this paper is to quantify and formalise these observations, which were suggested in \cite{j13error} relying partly on heuristic arguments.
Our asymptotic methods are based on the assumption that $\kappa$ and $\ee$ are small.
We calculate the probability that $|\x(t)|$ is less than $\ee$
in order to determine the influence of the smoothing function $F$ on the evolution of the system.
We also calculate the mean escape time of solutions from an $O(r)$ neighbourhood of the switching surface, where $\ee \ll r \ll 1$.
By restricting our analysis to a small neighbourhood we are able to reduce the necessary calculations to one dimension.


\section{Summary of mathematical results}
\label{sec:RESULTS}
\setcounter{equation}{0}

We are interested in what happens when a solution of (\ref{eq:genFilippov}) arrives at some point on the switching surface $\x=0$, when the perturbations of smoothing in (\ref{eq:genPerturbed}) and noise in (\ref{eq:genSDE}) are taken into account. Assuming that the solution has reached $\x=0$ in the perturbed model, and setting this point as the origin of coordinates, we study forward evolution of (\ref{eq:genSDE}) for positive times $t \le O(r)$, with
\begin{equation}
0 < \ee,\kappa \ll r \ll 1 \;.
\label{eq:eekappar}
\end{equation}
This implies $\bx(t) = O(r)$
because the drift of (\ref{eq:genSDE}) is bounded in an $O(r)$ neighbourhood.

It is helpful to treat sliding and crossing scenarios (\fig{fil}) separately.
Here we summarise our main results in each case.
The precise assumptions made for each result are stated in later sections.

The dynamics depends on the magnitude of the noise amplitude $\kappa$ relative to the size of the smoothing region $\ee$. If the origin belongs to a point on the switching surface where the Filippov convention predicts sliding, 
we calculate the probability $\mathbb{P}\left[ |\x(t)| \le \ee ~\big|~ \bx(0) = 0 \right]$ that a solution remains in the region $|\x| \le \ee$ after time $t$. We find:

\begin{enumerate}
\item
Suppose the origin belongs to an attracting sliding region of (\ref{eq:genFilippov})
and $\kappa\gg\sqrt\ee$.
Then for times in the intermediate range $O(\kappa^2) < t \le O(r)$, we have
\begin{equation}
\mathbb{P} \left[ |\x(t)| \le \ee ~\big|~ \bx(0) = 0 \right] = O \left( \frac{\ee}{\kappa^2} \right) \;.
\label{eq:Plarge}
\end{equation}

In this scenario, noise dominates over smoothing.
Solutions spend the majority of time evolving outside the neighbourhood of the switching surface, in $|\x(t)| \ge \ee$.
If $\lim_{\ee,\kappa \to 0} \frac{\kappa}{\sqrt{\ee}} = \infty$,
then, as with stochastically perturbed sliding motion \cite{simpson12},
in the zero noise limit we recover Filippov's sliding solution.
Note that we require $t > O(\kappa^2)$ so that the system has sufficient time
to settle to a quasi-steady-state distribution about $\x=0$.
We require $t \le O(r)$ to ensure $\bx(t) = O(r)$.

\item
Suppose the origin belongs to an attracting sliding region
and $\kappa\ll\sqrt\ee$.
Then for $t \le O(r)$ we have
\begin{equation}
\mathbb{P} \left[ |\x(t)| \le \ee ~\big|~ \bx(0) = 0 \right] = 1 - O \left( \frac{\kappa^2}{\ee} \right) \;.
\label{eq:Psmall}
\end{equation}

In this scenario, we find that evolution occurs primarily in the $\ee$ neighbourhood of the (smoothed out) switching surface, $|\x(t)| \le \ee$.
Therefore different choices of the smoothing function $F$ in (\ref{eq:genPerturbed}) lead to qualitatively different forward evolution. 
\end{enumerate}

The derivation of the above results are 
based on a steady-state approximation to the density of the $\x$-component of the solution about $\x=0$.
This approximation is ineffective or difficult to interpret
when the origin belongs to a crossing region, because
the piecewise-smooth vector field directs solutions away from $\x=0$.
For this reason, when Filippov's method predicts crossing we study the mean time for escape from $|\x| < r$, given by
\begin{equation}
\T = \mathbb{E} \left[ {\rm min} \left\{ t > 0 ~\big|~ |\x(t)| = r,\, \bx(0) = 0 \right\} \right] \;.
\label{eq:Tdef}
\end{equation}
We find:
\begin{enumerate}
\addtocounter{enumi}{2}
\item
If $\frac{\kappa}{\sqrt{\ee}} \ge O(1)$,
and the vector field points in the direction of increasing $\x$, then
\begin{equation}
\T = \frac{r}{a^+} + O(r^2) \;,
\label{eq:Tlarge}
\end{equation}
where $a^+ > 0$ is the $\x$-component of $f^+(0)$.
Note, in the absence of both smoothing and noise,
the forward orbit of the origin reaches $\x = r$ in a time $\frac{r}{a^+} + O(r^2)$.
Thus noise is sufficiently strong relative to smoothing to drive solutions quickly
out of the range of influence $|y|\le\ee$ of the smoothing function $F$.

\item
Again suppose the origin belongs to a crossing region
with the vector field pointing in the direction of increasing $\x$,
but now suppose $\kappa\ll\sqrt\ee$.
Here the nature of the smoothing function is important.
If $F \left( \frac{\x}{\ee}, 0 \right) > 0$ for all $\x \in [-\ee,\ee]$,
then we again have (\ref{eq:Tlarge}).
However, if $F \left( \frac{\x}{\ee}, 0 \right) < 0$ for an interval of $\x$ values, then
\begin{equation}
\T > O(r) \;.
\label{eq:Tsmall}
\end{equation}

Here, the vector field points to the left for part of $|\x| \le \ee$,
and the noise is highly unlikely
to drive solutions through this section and into $\x > \ee$ in an $O(r)$ time. Thus different smoothings of the piecewise-smooth system will lead to different dynamics, and may exhibit solutions that dwell near the discontinuity for large times at a crossing region, constituting spurious sliding dynamics. 
\end{enumerate}

It is not profitable to study $\T$ in a sliding scenario
because escape from $|\x| < r$ is extremely unlikely to occur over an $O(r)$ time frame,
regardless of the relative size of $\ee$ and $\kappa$ and the nature of the smoothing function. 

\section{Reduction to one dimension and scaling}
\label{sec:SCALE}
\setcounter{equation}{0}

Under assumptions that are stated below, the dynamics of (\ref{eq:genSDE}) in the $\x$-direction is described by a one-dimensional stochastic differential equation.

As explained in \S\ref{sec:RESULTS}, we consider 
$\bx(t) = O(r)$, and therefore we can write (\ref{eq:genSDE}) as
\begin{equation}
d\bx(t) = 
\left\{ \begin{array}{lc}
f^-(0) \;, & \x(t) \le -\ee \\
F \left( \frac{\x(t)}{\ee}, 0 \right)  \;, & -\ee < \x(t) < \ee \\
f^+(0) \;, & \x(t) \ge \ee
\end{array} \right\} \,dt + \kappa D(0) \,d\bW(t)+ O(r) \;.
\label{eq:scaledSDE}
\end{equation}
Since the right-hand side of (\ref{eq:scaledSDE}) depends only on $\frac{\x}{\ee}$ to leading order,
we are able to perform an analysis of the stochastic dynamics in $\x$ independently of the remaining components of $\bx$.
The $y$-component of (\ref{eq:scaledSDE}) is
\begin{equation}
d\x(t) = 
\left\{ \begin{array}{lc}
a^-  \;, & \x(t) \le -\ee \\
A \left( \frac{\x(t)}{\ee} \right)  \;, & -\ee < \x(t) < \ee \\
a^+  \;, & \x(t) \ge \ee
\end{array} \right\} \,dt + \kappa  B^T \,d\bW(t)+ O(r) \;,
\label{eq:scaledSDEx}
\end{equation}
where $a^{\pm}$ is the $\x$-component of $f^{\pm}(0)$,
$A\left(\frac{\x}{\ee}\right)$ is the $\x$-component of $F\left(\frac{\x}{\ee},0\right)$,
and $B^{\sf T}$ is the appropriate row of $D(0)$.

The noise in (\ref{eq:scaledSDEx}) is a sum of $n$ independent Brownian motions, equivalent to a single Brownian motion of amplitude $\sqrt{|B|} + O(r)$.
Let us absorb $\sqrt{|B|}$ into the value of $\kappa$, and define
\begin{equation}
\phi\left(\frac{\x}{\ee}\right) = \left\{ \begin{array}{lc}
a^- \;, & \x\le -\ee \\
A\left(\frac{\x}{\ee}\right) \;, & -\ee < \x< \ee \\
a^+ \;, & \x \ge \ee
\end{array} \right. \;.
\label{eq:phi}
\end{equation}
Then, neglecting $O(r)$ terms, (\ref{eq:scaledSDEx}) becomes
\begin{equation}
d\x(t) = \phi \left( \frac{\x}{\ee} \right) \,dt + \kappa \,dW(t) \;.
\label{eq:reducedSDE}
\end{equation}
The smoothing of the system over an $O(\ee)$ spacial and temporal scale suggests the scaling
\begin{equation}
\tilde{\x} = \frac{\x}{\ee} \;, \qquad
\tilde{t} = \frac{t}{\ee} \;, \qquad
\tilde{\kappa} = \frac{\kappa}{\sqrt{\ee}} \;,
\label{eq:tildexykappa}
\end{equation}
with which (\ref{eq:reducedSDE}) becomes
\begin{equation}
d\tilde{\x}(\tilde{t}) = \phi(\tilde{\x}(\tilde{t})) \,d\tilde{t} + \tilde{\kappa} \,dW(\tilde{t}) \;.
\label{eq:reducedSDEtilde}
\end{equation}
Equation (\ref{eq:reducedSDEtilde}) is the reduced system that we analyse in the next two sections.
The signs of $a^-$ and $a^+$ determine whether dynamics at the origin of the unperturbed system (\ref{eq:genFilippov}) involves sliding (or attracting or repelling type) or crossing.
Ignoring tangencies, we have the classification:
\begin{center}
\begin{tabular}{|l|c|c|}\hline
& $a^+ < 0$ & $a^+ > 0$ \\
\hline
$a^- < 0$ & crossing & sliding (repelling) \\
\hline
$a^- > 0$ & sliding (attracting) & crossing \\\hline
\end{tabular}
\end{center}
It is sufficient for our purposes to consider only the lower row of this table.

\section{Steady-state approximations and the mean escape time}
\label{sec:STOCH}
\setcounter{equation}{0}

\subsection{Steady-state approximations}
\label{sub:SS}

\begin{figure}[h!]
\begin{center}
\includegraphics[width=\textwidth]{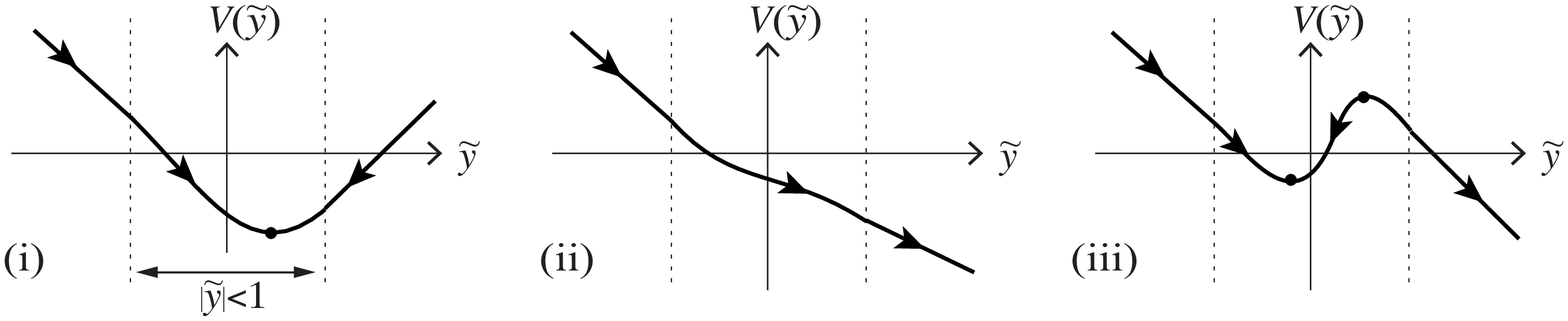}
\caption{\sf\small Sketches of the potential $V(\tilde \x)$ corresponding to the vector fields in \fig{flows}, showing:
(i) a global potential well, (ii) no potential well, (iii) a local potential well.
\label{fig:potentialSketch}
}
\end{center}
\end{figure}

The Fokker-Planck equation for (\ref{eq:reducedSDEtilde}) (see e.g. \cite{schuss10,g09}) is
\begin{equation}
\frac{\partial p(\tilde{\x},\tilde{t})}{\partial \tilde{t}} = \frac{\partial}{\partial \tilde{\x}}
\left( -\phi(\tilde{\x}) p(\tilde{\x},\tilde{t}) + \frac{\tilde{\kappa}^2}{2} \frac{\partial p(\tilde{\x},\tilde{t})}{\partial \tilde{\x}} \right) \;,
\label{eq:FPE}
\end{equation}
where $p(\tilde\x,\tilde t)$ denotes the transitional probability density function of (\ref{eq:reducedSDEtilde}).
For a steady-state density,
$p(\tilde{\x},\tilde{t}) = p_{\rm ss}(\tilde{\x})$, the left-hand side of (\ref{eq:FPE}) is zero.
Assuming $p_{\rm ss}$ vanishes as $\tilde{\x} \to \pm \infty$, if $p_{\rm ss}$ exists it must be given by
\begin{equation}
p_{\rm ss}(\tilde{\x}) = K {\rm e}^{\frac{-2 V(\tilde{\x})}{\tilde{\kappa}^2}} \;,
\label{eq:pss}
\end{equation}
where
\begin{equation}
V(\tilde{\x}) = -\int_{-1}^{\tilde{\x}} \phi(v) \,dv \;,
\label{eq:V}
\end{equation}
represents a potential function, and $K$ is a normalisation constant.
By (\ref{eq:phi}) we can write
\begin{equation}
V(\tilde{\x}) = \left\{ \begin{array}{lc}
-a^-(\tilde{\x}+1) \;, & \tilde{\x} \le -1 \\
-\int_{-1}^{\tilde{\x}} \phi(v) \,dv \;, & -1 \le \tilde{\x} \le 1 \\
V(1) - a^+(\tilde{\x}-1) \;, & \tilde{\x} \ge 1
\end{array} \right. \;.
\label{eq:V2}
\end{equation}
Consequently (\ref{eq:pss}) is normalisable only if $a^- > 0$ and $a^+ < 0$, corresponding to sliding as in \fig{potentialSketch}(i).
The reduced system (\ref{eq:reducedSDEtilde}) has a steady-state density
exactly when the origin corresponds to a attracting sliding region of (\ref{eq:genFilippov}).

In all other cases (i.e.~crossing or repelling sliding),
(\ref{eq:reducedSDEtilde}) does not have a steady-state density
because for $|\tilde{\x}| > 1$ the drift directs solutions away from $\tilde\x=0$ on at least one side.
However, $V$ may have a local potential well,
as in \fig{potentialSketch}(iii),
or perhaps many potential wells.
For \fig{potentialSketch}(iii), if the noise amplitude is sufficiently small,
then with high probability solutions become trapped in the well for a relatively long period of time.
For an initial condition in the well,
after a brief transient phase, the probability density function $p$ takes a near-steady-state form
for which probability leaks out of the well on a long time-scale.
In this scenario it is more useful to consider the mean escape time from the region $|\tilde{\x}| \le 1$. 


\subsection{Mean escape time}
\label{sub:MET}

We let $\tilde{\T}$ denote the mean escape time for (\ref{eq:reducedSDEtilde}), defined as
\begin{equation}
\tilde{\T} = \mathbb{E} \left[ {\rm min} \left\{ \tilde{t} > 0 ~\big|~
|\tilde{\x}(\tilde{t})| = \tilde{r},\, \tilde{\x}(0) = 0 \right\} \right] \;,
\label{eq:Ttildedef}
\end{equation}
where 
\begin{equation}
\tilde{r} = \frac{r}{\ee} \;.
\label{eq:tilder}
\end{equation}
To find the mean escape time to $\tilde{\x} = \tilde{r}$,
we consider the Fokker-Planck equation (\ref{eq:FPE}),
together with the initial and boundary conditions
\begin{equation}
p(\tilde{\x},0) = \delta(\tilde{\x}) \;, \qquad
p(\pm \tilde{r},\tilde{t}) = 0 \;.
\label{eq:BCs2}
\end{equation}
The first condition in (\ref{eq:BCs2}) simply corresponds to fixing $\tilde{\x}(\tilde{t}) = 0$.
The second condition in (\ref{eq:BCs2}) represents absorption at $\tilde{\x} = \tilde{r}$.
It is a standard exercise to demonstrate that the boundary value problem 
(\ref{eq:FPE}) with (\ref{eq:BCs2}), implies that the mean escape time is given by
\begin{equation}
\tilde{\T} = \int_0^\infty \int_{-\tilde{r}}^{\tilde{r}} p(\tilde{\x},\tilde{t}) \,d\tilde{\x} \,d\tilde{t} \;,
\label{eq:Ttilde0}
\end{equation}
see for instance \cite{schuss10,g09}.
Since (\ref{eq:reducedSDEtilde}) is one-dimensional and time-independent,
we can obtain an explicit expression for $\tilde{\T}$.
By integrating (\ref{eq:FPE}) over all positive time
and solving the resulting ordinary differential equation, 
we obtain
\begin{equation}
\tilde{\T} = \frac{2}{\tilde{\kappa}^2} \int_{-\tilde{r}}^{\tilde{r}} \int_{-\tilde{r}}^v
\left( H(v) - C \right)
{\rm e}^{\frac{2(V(v)-V(u))}{\tilde{\kappa}^2}}
\,du \,dv \;,
\label{eq:Ttilde}
\end{equation}
where $H(v)$ is the Heaviside function, and $C$ is given by
\begin{equation}
C = \frac{\int_0^{\tilde{r}} {\rm e}^{\frac{2 V(w)}{\tilde{\kappa}^2}} \,dw}
{\int_{-\tilde{r}}^{\tilde{r}} {\rm e}^{\frac{2 V(w)}{\tilde{\kappa}^2}} \,dw} \;.
\label{eq:C}
\end{equation}
Formulae similar to (\ref{eq:Ttilde}) are derived in \cite{schuss10,g09,gh99,s51}.
For convenience we provide a derivation of (\ref{eq:Ttilde}) in Appendix \ref{sec:TDERIV}.

\section{Asymptotics for the occupation probability and the mean escape time}
\label{sec:CALC}
\setcounter{equation}{0}

\subsection{Sliding}
\label{sec:STABLE}
\setcounter{equation}{0}

Here we suppose $a^- > 0$ and $a^+ < 0$ and investigate the steady-state probability density function $p_{\rm ss}$ given by (\ref{eq:pss}).
If $\tilde{\x}(\tilde{t})$ is distributed according to (\ref{eq:pss}), then
\begin{equation}
\mathbb{P} \left[ |\tilde{\x}(\tilde{t})| \le 1 \right] =
\frac{\int_{-1}^1 {\rm e}^{\frac{-2 V(u)}{\tilde{\kappa}^2}} \,du}
{\int_{-\infty}^\infty {\rm e}^{\frac{-2 V(u)}{\tilde{\kappa}^2}} \,du} \;.
\label{eq:P}
\end{equation}
From (\ref{eq:V2}) we determine
\begin{equation}
\int_{-\infty}^{-1} {\rm e}^{\frac{-2 V(u)}{\tilde{\kappa}^2}} \,du =
\frac{\tilde{\kappa}^2}{2 a^-} \;, \qquad
\int_1^{\infty} {\rm e}^{\frac{-2 V(u)}{\tilde{\kappa}^2}} \,du =
\frac{\tilde{\kappa}^2}{-2 a^+} \,{\rm e}^{\frac{-2 V(1)}{\tilde{\kappa}^2}} \;.
\end{equation}
Dividing (\ref{eq:P}) through by its numerator gives
\begin{equation}
\mathbb{P} \left[ |\tilde{\x}(\tilde{t})| \le 1 \right] =
\frac{1}{1 + \frac{\frac{\tilde{\kappa}^2}{2 a^-} + \frac{\tilde{\kappa}^2}{-2 a^+} \,{\rm e}^{\frac{-2 V(1)}{\tilde{\kappa}^2}}}
{\int_{-1}^1 {\rm e}^{\frac{-2 V(u)}{\tilde{\kappa}^2}} \,du}} \;.
\label{eq:P2}
\end{equation}
In the case $\tilde{\kappa} \gg 1$, we can then simply use
${\rm e}^{-\frac{2 V(u)}{\tilde{\kappa}^2}} = 1 + O \left( \frac{1}{\tilde{\kappa}^2} \right)$,
to find
\begin{equation}
\mathbb{P} \left[ |\tilde{\x}(\tilde{t})| \le 1 \right] =
\frac{4}{\left( \frac{1}{a^-} + \frac{1}{-a^+} \right) \tilde{\kappa}^2} +
O \left( \frac{1}{\tilde{\kappa}^4} \right) \;.
\label{eq:Plarge2}
\end{equation}
The case $\tilde{\kappa} \ll 1$ requires a more detailed argument to obtain a bound on $\mathbb{P}$.
If $V(1) \le 0$ (the case $V(1) > 0$ may be treated analogously),
since $V$ is differentiable with $V'(1) = -a^+ > 0$,
there exists a value $0 < L < 2$ such that $V(u) \le V(1)$ for all $1-L \le u \le 1$.
Consequently 
$\int_{-1}^1 {\rm e}^{\frac{-2 V(u)}{\tilde{\kappa}^2}} \,du \ge
L {\rm e}^{\frac{-2 V(1)}{\tilde{\kappa}^2}}$.
Then from (\ref{eq:P2}),
\begin{equation}
\mathbb{P} \left[ |\tilde{\x}(\tilde{t})| \le 1 \right] \ge
\frac{1}{1 + \frac{\tilde{\kappa}^2}{2 L}
\left( \frac{1}{a^-} \,{\rm e}^{\frac{2 V(1)}{\tilde{\kappa}^2}} +
\frac{1}{-a^+} \right)} \;.
\end{equation}
Therefore for $\tilde{\kappa} \ll 1$ we have
\begin{equation}
\mathbb{P} \left[ |\tilde{\x}(\tilde{t})| \le 1 \right] = 1 + O(\tilde{\kappa}^2) \;.
\label{eq:Psmall2}
\end{equation}

We now interpret (\ref{eq:Plarge2}) and (\ref{eq:Psmall2}) for the original system (\ref{eq:genSDE}).
An $O(r)$ error in the coefficients of (\ref{eq:reducedSDEtilde})
relates to an $O(r)$ error in the potential $V$.
Consequently an $O(r)$ error does not contribute additively to $\mathbb{P}$,
rather it appears in the coefficients of the terms in (\ref{eq:Plarge2}) and (\ref{eq:Psmall2}).
Hence this error appears as higher order contributions in both cases,
and we have (\ref{eq:Plarge}) and (\ref{eq:Psmall}) of \S\ref{sec:RESULTS}.
A formal demonstration of this argument is beyond the scope of this paper.

For the case $\tilde{\kappa} \gg 1$, in view of (\ref{eq:Plarge2}),
the system (\ref{eq:reducedSDEtilde}) behaves like Brownian motion
of amplitude $\tilde{\kappa}$ with piecewise-constant drift,
for which the correlation time is $O(\tilde{\kappa}^2)$ \cite{simpson12}.
Therefore (\ref{eq:reducedSDEtilde}) settles to
the steady-state density on an $O(\tilde{\kappa}^2)$ time-scale,
and for this reason we assume $t > O(\kappa^2)$ for the result (\ref{eq:Plarge}).
No lower bound on $t$ is given for (\ref{eq:Psmall}), corresponding to $\tilde{\kappa} \ll 1$,
because in this case for very short times, although (\ref{eq:reducedSDEtilde}) will
have not settled to steady-state, solutions will reside inside the smoothing region $|\tilde\x|<1$ with high probability.

\subsection{Crossing}
\label{sub:CROSSING}

Here we suppose $a^-, a^+ > 0$, and analyse the asymptotics of $\tilde{\T}$ as given by (\ref{eq:Ttilde}).
When $a^-, a^+ > 0$, straight-forward bounding arguments applied to (\ref{eq:C}) (see Appendix \ref{sec:C} for details)
reveal that $C$ is exponentially small
(a simple upper bound is $C \le 2 a^- \left( \frac{1}{a^+} + \frac{2}{\tilde{\kappa}^2} \right)
{\rm e}^{\frac{- \tilde{r}a^-}{\tilde{\kappa}^2}}$).
This is because terms involving $C$ relate to escape through $\tilde{\x} = -1$, which is highly unlikely.
Indeed, terms involving $C$ constitute high order contributions in the manipulations that follow, and do not appear in our final results.
Hence for simplicity here we take $C = 0$,
with which (\ref{eq:Ttilde}) is given by
\begin{equation}
\tilde{\T} = \frac{2}{\tilde{\kappa}^2} \int_0^{\tilde{r}}
{\rm e}^{\frac{2 V(v)}{\tilde{\kappa}^2}} \int_{-\tilde{r}}^v
{\rm e}^{\frac{-2 V(u)}{\tilde{\kappa}^2}}
\,du \,dv \;.
\label{eq:TtildeCzero}
\end{equation}
To evaluate (\ref{eq:TtildeCzero}) asymptotically,
we split the two-dimensional domain of integration into four pieces:
\begin{align}
\tilde{\T} &=
\frac{2}{\tilde{\kappa}^2}\left\{\;\; 
\int_0^1  \int_{-\tilde{r}}^{-1} 
+ \int_0^1 \int_{-1}^v 
+ \int_1^{\tilde{r}}  \int_{-\tilde{r}}^1
+ \int_1^{\tilde{r}}  \int_1^v  \;\;\right\} 
{\rm e}^{\frac{2\left(V(v)- V(u)\right)}{\tilde{\kappa}^2}} \,du \,dv \;.
\label{eq:fourPieces}
\end{align}
If $\tilde{\kappa} \ge O(1)$, i.e.~$\tilde{\kappa}$ is not small,
via explicit integration and applying simple bounds
we find that the first three double integrals of (\ref{eq:fourPieces}) are at most $O(\tilde{\kappa}^2)$,
and because $\tilde{r} \gg \tilde{\kappa}^2$ from (\ref{eq:eekappar}), these
are dominated by the fourth double integral:
\begin{equation}
\frac{2}{\tilde{\kappa}^2} \int_1^{\tilde{r}} {\rm e}^{\frac{2 V(v)}{\tilde{\kappa}^2}} \int_1^v {\rm e}^{\frac{-2 V(u)}{\tilde{\kappa}^2}} \,du \,dv
= \frac{2}{\tilde{\kappa}^2} \int_1^{\tilde{r}} {\rm e}^{\frac{-2 a^+ v}{\tilde{\kappa}^2}} \int_1^v {\rm e}^{\frac{2 a^+ u}{\tilde{\kappa}^2}} \,du \,dv
= \frac{\tilde{r}}{a^+} + O(1,\tilde{\kappa}^2) \;,
\label{eq:fourthTerm}
\end{equation}
Hence, if $\tilde{\kappa} \ge O(1)$,
\begin{equation}
\tilde{\T} = \frac{\tilde{r}}{a^+} \left( 1 + O \left( \frac{\kappa^2}{r} \right) \right) \;.
\label{eq:Ttildelarge}
\end{equation}

Alternatively, if $\tilde{\kappa}$ is small
then $A(\tilde{\x})$, which appears in (\ref{eq:phi}), is important.
Recall that $A(\tilde{\x})$ is continuous for $\tilde{\x} \in [-1,1]$ and that $A(\pm 1) = a^\pm > 0$.
If $A(\tilde{\x})$ is positive throughout the interval $[-1,1]$,
then $\tilde{\T} \to \frac{\tilde{r}}{a^+}$ as $\ee \to 0$ and $\kappa \to 0$, as one would expect.
Let us suppose that $A(\tilde{\x}) < 0$ for some part of the interval $[-1,1]$.
Furthermore, suppose $A(\tilde{\x})$ has exactly two roots, $-1 < \tilde{\x}_1 < \tilde{\x}_2 < 1$,
and $\frac{\partial A(\tilde{\x}_1)}{\partial \tilde{\x}} < 0$ and
$\frac{\partial A(\tilde{\x}_2)}{\partial \tilde{\x}} > 0$.
The potential $V(\tilde{\x})$ then has a well (as in \fig{potentialSketch}(iii)),
with a local minimum at $\tilde{\x}_1$ and a local maximum at $\tilde{\x}_2$.
For simplicity suppose $\tilde{\x}_2 > 0$,
so that solutions must pass through the well in order to reach $\tilde{\x} = \tilde{r}$.
This is sensible because for crossing we expect solutions to originate from negative $\x$, at least locally.

We can evaluate (\ref{eq:TtildeCzero}) asymptotically 
by treating $\tilde{\kappa}$ and $\frac{1}{\tilde{r}}$ as independent small parameters.
The result depends on the relative size of
$\tilde{\kappa}$ and $\frac{1}{\tilde{r}}$ as the limit is taken.
For instance, if we take $\tilde{\kappa} \to 0$,
the maximum contribution to the integral comes from the point at which $V(v)-V(u)$ attains its maximum
(the point $(u,v) = (\tilde{\x}_1,\tilde{\x}_2)$).
Alternatively, if we take $\frac{1}{\tilde{r}} \to 0$,
the size of the domain of integration tends to infinity and the maximum contribution relates to
$(u,v)$ far from $(0,0)$.

Let us consider (\ref{eq:fourPieces}) and evaluate the four double integrals separately, assuming as above that $A(\tilde{\x})$ has two zeros
$$A(\tilde{\x_{1,2}})=0\quad:\quad -1 < \tilde{\x}_1 < \tilde{\x}_2 < 1,\quad\frac{\partial A(\tilde{\x}_1)}{\partial \tilde{\x}} < 0<\frac{\partial A(\tilde{\x}_2)}{\partial \tilde{\x}} \;.$$
The point $(u,v) = (\tilde{\x}_1,\tilde{\x}_2)$ lies in the domain of integration
of the second integral.
A straight-forward application of Laplace's method (see Appendix \ref{sec:doubleint}) yields
\begin{equation}
\frac{2}{\tilde{\kappa}^2} \int_0^1 {\rm e}^{\frac{2 V(v)}{\tilde{\kappa}^2}}
\int_{-1}^v {\rm e}^{\frac{-2 V(u)}{\tilde{\kappa}^2}} \,du \,dv
= \frac{2 \pi}
{\sqrt{-A'(\tilde{\x}_1) A'(\tilde{\x}_2)}}
\,{\rm e}^{\frac{2}{\tilde{\kappa}^2} \left( V(\tilde{\x}_2) - V(\tilde{\x}_1) \right)}
\big( 1 + O(\tilde{\kappa}) \big) \;.
\label{eq:secondTerm}
\end{equation}
Repeating (\ref{eq:fourthTerm}),
we find that the fourth double integral is given by
\begin{equation}
\frac{2}{\tilde{\kappa}^2} \int_1^{\tilde{r}} {\rm e}^{\frac{2 V(v)}{\tilde{\kappa}^2}} \int_1^v {\rm e}^{\frac{-2 V(u)}{\tilde{\kappa}^2}} \,du \,dv
= \frac{\tilde{r}}{a^+} + O(1) \;,
\label{eq:fourthTerm2}
\end{equation}
where here the $O(\tilde{\kappa}^2)$ term in (\ref{eq:fourthTerm}) is of higher order than the error term given in (\ref{eq:fourthTerm2}).
The remaining two double integrals of (\ref{eq:fourPieces}) only provide higher order contributions
that may be absorbed into the error term of (\ref{eq:secondTerm}).
This can be verified via a combination of explicit integration and Laplace's method.

In summary, when $\tilde{\kappa} \ll 1$, the term (\ref{eq:secondTerm}) appears because there is a potential well.
This term corresponds to a large increase in the mean escape time due to solutions becoming trapped in the well.
In contrast, the term (\ref{eq:fourthTerm2}) corresponds to the time taken by purely following the vector field $f^+$,
and arises independently to the nature of the smoothing function.
We can therefore write
\begin{equation}
\tilde{\T} = \frac{2 \pi {\cal S}}
{\sqrt{-A'(\tilde{\x}_1) A'(\tilde{\x}_2)}}
\,{\rm e}^{\frac{2}{\tilde{\kappa}^2} \left( V(\tilde{\x}_2) - V(\tilde{\x}_1) \right)} \big( 1 + O(\tilde{\kappa}) \big)
+ \frac{\tilde{r}}{a^+} \left( 1 + O \left(\mbox{$ \frac{1}{\tilde{r}} $} \right) \right) \;,
\label{eq:Ttildesmall}
\end{equation}
where ${\cal S}$ is a {\it Stokes multiplier} (see e.g. \cite{heading,bo99,berry-uniform}) whose value is
\begin{equation}
{\cal S}=\left\{\begin{array}{lll}1&&\mbox{if $V$ has a potential well}
\;,\\0&&\mbox{if $V$ has no potential well}\;. \end{array}\right.
\label{eq:S}
\end{equation}
For $V$ to have a potential well, the values of $\tilde \x_1$ and $\tilde \x_2$ must be real and lie inside the neighbourhood of the switching surface given by $-1<\tilde \x_1<1$ and $0<\tilde \x_2<1$. 

The last step is to interpret (\ref{eq:Ttildelarge}) and (\ref{eq:Ttildesmall}) for the original system (\ref{eq:genSDE}).
For the case $\tilde{\kappa} \gg 1$, as in the previous section,
we note that the $O(r)$ error in the coefficients of (\ref{eq:reducedSDEtilde})
relate to an $O(r)$ error in the potential $V$.
In (\ref{eq:fourthTerm}) this leads to an $O(r)$ error in the coefficient of the leading order term.
That is,
\begin{equation}
\T = \ee \left( 1 + O(r) \right) \tilde{\T} \;,
\label{eq:TtildetoT}
\end{equation}
where the $\ee$ is present because $t = \ee \tilde{t}$, and we therefore have (\ref{eq:Tlarge}).
For the case $\tilde{\kappa} \ll 1$, equation (\ref{eq:Ttildesmall}) indicates that $\tilde{\T}$ is exponentially large.
Since the reduced system (\ref{eq:reducedSDEtilde}) only applies for $t \le O(r)$, however, 
we can conclude only that $\T > O(r)$, giving the statement (\ref{eq:Tsmall}).
As in the previous subsection, more detailed results are beyond the scope of this paper.

\section{Example: a dry-friction oscillator}
\label{sec:FRICTION}
\setcounter{equation}{0}

\subsection{Modelling assumptions}

\begin{figure}[h!]
\begin{center}
\setlength{\unitlength}{1cm}
{\includegraphics[width=0.5\textwidth]{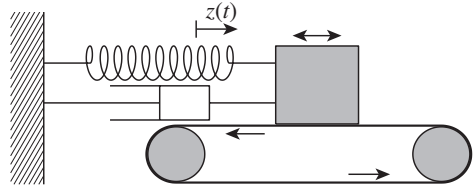}}
\caption{\sf\small The friction oscillator. A block is subject to linear spring and damping forces, and sits on a belt moving at constant unit speed to the left. The spring extension is $z(t)$. 
\label{fig:frictionOsc}
}
\end{center}
\end{figure}

To demonstrate the effect of the results above, let us consider a toy model of a dry-friction oscillator \begin{equation}
\ddot{z} + \dot{z} + z + \cF(\dot{z}+1) =0 \;,
\label{eq:osc1}
\end{equation}
as shown in \fig{frictionOsc}, 
where $z$ denotes the time-dependent horizontal displacement of a block of unit mass,
relative to a belt moving at unit speed to the left.
The block is subject to a damping coefficient and a spring constant which are both unity,
and to a velocity-dependent dry-friction force $\cF$ 
due to the contact of the block with the belt.
The oscillator is a prototypical model
useful for studying frction dynamics in diverse mechanical systems \cite{popp96,fm94}.

\begin{figure}[h!]
\begin{center}
\includegraphics[width=0.5\textwidth]{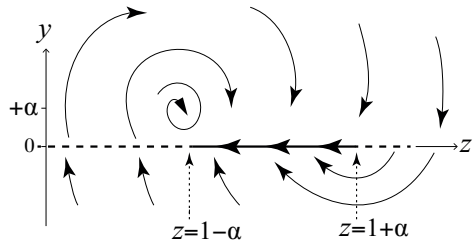}
\caption{\sf\small 
A sketch of the phase portrait of (\ref{eq:osc2}), as described in the text, showing dynamics of the friction oscillator according to Filippov's method. There are two crossing regions (dashed), around a sliding region (solid line) where the sliding vector field points to the left as shown. Dynamics in the regions $\x>0$, $\x<0$, and $\x=0$, correspond mechanically to leftward slip, rightward slip, and sticking.
\label{fig:phase}
}
\end{center}
\end{figure}

We let $\sp = \dot{z} + 1$ denote the velocity of the block relative to the belt,
with which (\ref{eq:osc1}) may be written as the first order system
\begin{equation}
\begin{split}
\dot{z} &= \sp - 1 \;, \\
\dot{\sp} &= 1 - z - \sp - \cF(\sp) \;.
\end{split}
\label{eq:osc2}
\end{equation}
The simplest Coulomb model of the friction force is
\begin{equation}
\cF(\sp) = \alpha \,{\rm sgn}(\sp) \;,
\label{eq:forcePWC}
\end{equation}
for some constant $\alpha > 0$. 
For (\ref{eq:osc2}) with (\ref{eq:forcePWC}), $\sp=0$ is a switching surface
which, in the standard Filippov convention, is an attracting sliding region for $1-\alpha < z < 1+\alpha$,
and a crossing region for $z < 1-\alpha$ and $z > 1+\alpha$,
see Fig.~\ref{fig:phase}. Sliding corresponds to a mechanical sticking phase, in which the relative speed remains at $\x=0$ for an interval of time, and crossing corresponds to instantaneous switch between slip-to-the-left and slip-to-the-right. 

\begin{figure}[h!]
\begin{center}
\includegraphics[width=0.5\textwidth]{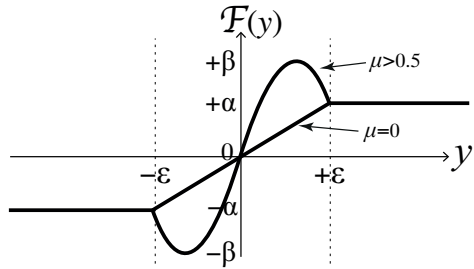}
\caption{\sf\small 
A sketch of the friction force ${\cal F}(y)$. For $\mu=0$ $\cal F$ is linear in $|\sp|<\ee$. For $\mu>1/2$ there exist two turning points in $|\sp| < \ee$ at $\sp_s=\pm\ee\sqrt{(1+\mu)/3\mu}$. 
\label{fig:fricforce}
}
\end{center}
\end{figure}

The friction model (\ref{eq:forcePWC}) neglects any dependence of the contact force $\cal F$ on the speed, i.e. the absolute value of $y$. Our only concern here is the value for $y\approx0$, specifically the difference between the kinetic friction force ${\cal F}=\pm\alpha$ when ${\rm sgn}(y)=\pm1$, and the static friction force when $y=0$, at which the value of (\ref{eq:forcePWC}) is ambiguous. To resolve this, let us consider a well-defined perturbation of (\ref{eq:forcePWC}) where
\begin{equation}
\cF(\sp) =\alpha
\left\{ \begin{array}{lc}
-1 \;, & \sp \le -\ee \\
\frac{\sp}{\ee} + \mu (\frac{\sp}{\ee}-(\frac{\sp}{\ee})^3) \;, & -\ee < \sp < \ee \\
1 \;, & \sp \ge \ee
\end{array} \right. \;,
\label{eq:forcePerturbed}
\end{equation}
with $\alpha > 0$ and $\mu \ge 0$. This provides a continuous change in the friction force across $\sp=0$
by spreading the jump out across a region $|\sp|\le\ee$ for small $\ee$.
For $\mu>\frac{1}{2}$, $\cF$ has turning points 
at $\sp_s=\pm\ee\sqrt{\frac{1+\mu}{3\mu}}$, (i.e.~$\cF'(\sp_s)=0$).
We also let
\begin{equation}
\beta \equiv \max_{|\sp|\le\ee} \,[ \cF(\sp) ] =
\alpha \left\{ \begin{array}{lll}
\frac{2 (1+\mu)^{3/2}}{3\sqrt{3\mu}} & {\rm for} & \mu>1/2 \;,\\
1 & {\rm for} & \mu<1/2\;,
\end{array}\right.
\end{equation}
and note that $\frac{2(1+\mu)^{3/2}}{3\sqrt{3\mu}}>1$ for $\mu>\frac12$.
A sketch of $\cF$ is shown in \fig{fricforce}.
For $\mu > \frac{1}{2}$ the perturbation (\ref{eq:forcePerturbed}) introduces a breakaway force $\pm\beta$ that exceeds the force $\pm\alpha$ required to keep the object in motion.

Let us now consider the effect of random errors in the model, represented by white noise, $\kappa \frac{dW}{dt}$, added to the force on the block, and assume that the noise amplitude $\kappa$ is small.
The oscillator is then described by the following system of stochastic differential equations:
\begin{equation}
\begin{split}
dz(t) &= (\sp(t)-1) \,dt \;, \\
 \,d\sp(t) &= \big( 1 - z(t) - \sp(t) - \cF(\sp) \big) \,dt + \kappa \,dW(t) \;.
\end{split}
\label{eq:osc3}
\end{equation}
With $\bx = (\sp,z)$, the system (\ref{eq:osc3}) is of the form (\ref{eq:genSDE}).

\subsection{Calculations for transitions from stick to slip}

To consider the escape of the system from a neighbourhood of $\sp = 0$, 
let us consider an initial state $(\sp,z) = (0,z_0)$.

For (\ref{eq:osc3}),
the leading order approximation to dynamics orthogonal to the switching surface,
derived in general in \S\ref{sec:SCALE}, and valid over short times, is given by
\begin{equation}
d\sp(t) = \phi \left( \frac{\sp}{\ee} \right) \,dt + \kappa \,dW(t) \;,
\end{equation}
where
\begin{equation}
\phi(u) = \left\{ \begin{array}{lc}
1-z_0+\alpha \;, & u\le -1 \\
1-z_0-\alpha(u+\mu(u-u^3))  \;, & -1 < u< 1 \\
1-z_0-\alpha \;, & u \ge 1
\end{array} \right. \;.
\label{eq:phifric0}
\end{equation}
Integrating to find the potential (\ref{eq:V}) gives
\begin{equation}
V\left(u\right) = 
-\left\{ \begin{array}{lc}
(1-z_0+\alpha)(u+1) \;, & u\le -1 \\
(1-z_0)(u+1) - \frac{\alpha (1+\mu)(u^2-1)}{2} + \frac{\alpha \mu (u^4-1)}{4}  \;, & -1 < u< 1 \\
2(1-z_0)+(1-z_0-\alpha)(u-1) \;, & u \ge 1
\end{array} \right. \;.
\label{eq:Vfric}
\end{equation}

Below we use the results of \S\ref{sub:MET} and \S\ref{sub:CROSSING}
to determine the mean time for escape from a neighbourhood of $\sp = 0$.
We do not calculate the probability density function in the sliding region, but remark that, for this model,
it coincides with Filippov's convention in the zero noise limit.


We instead consider the region $z_0<1-\alpha$, where Filippov's convention predicts crossing. 
The potential $V$ has turning points where $\phi$ vanishes, namely at $u=u_i$, defined by
\begin{equation}\label{cubic}
\phi(u_i)=1-z_0-\alpha(u_i+\mu (u_i- u_i^3))=0\;\;\;{\rm and}\;\;\; |u_i|<1\;,
\end{equation}
which can have up to three real solutions $u_1,u_2,u_3$. 
Two turning points coalesce at
\begin{equation}		
u=u^{(\pm)}:=\pm\sqrt{\frac{1+\mu
}{3\mu}},\qquad z_0=z_0^{(\pm)}:=1\mp\frac{2\alpha(1+\mu
)^{3/2}}{3\sqrt{3\mu}}\;,
\end{equation}
which lie inside $|u|<1$ for $\mu>1/2$, and which correspond to the maxima and minima of $\cal F$ in (\ref{eq:forcePerturbed}), noting $z_0^{(\pm)}=1\mp\beta$. 
The point $z=z_0^{(+)}$ forms a breakaway point at which the block is released from sticking into rightward slip, i.e. from $u=0$ into $u>0$. 

For $z_0^{(+)}\le z<1-\alpha$, and parameter values of physical interest, there are two roots of (\ref{cubic}), $u_1$ and $u_2$, that we label such that $-1<u_1\le u_2<1$,
and for which $\phi'(u_1) \le 0 \le \phi'(u_2)$.
The function $V(v)-V(u)$, about which the second term of (\ref{eq:Ttildesmall}) is approximated,
has a maximum at $\left\{u,v\right\}=\left\{u_1,u_2\right\}$.
Specifically, in (\ref{eq:Ttildesmall}) we have
\begin{equation}\label{fr1}
\begin{array}{rcl}
a^+&=&1-z_0+\alpha \;,\\
A'(\tilde \x_i) &=& -\alpha(1+\mu-3\mu\tilde \x_i^2) \;,\quad i=1,2,
\end{array}
\end{equation}
where $\tilde{\x}_i=u_i$, and
\begin{equation}\label{fr2}
\phi(u_i)=0 \;, \qquad
\phi'(u_1) \le 0 \le \phi'(u_2) \;.
\end{equation}

\subsection{Numerical computations}

\begin{figure}[h!]
\begin{center}
\includegraphics[width=\textwidth]{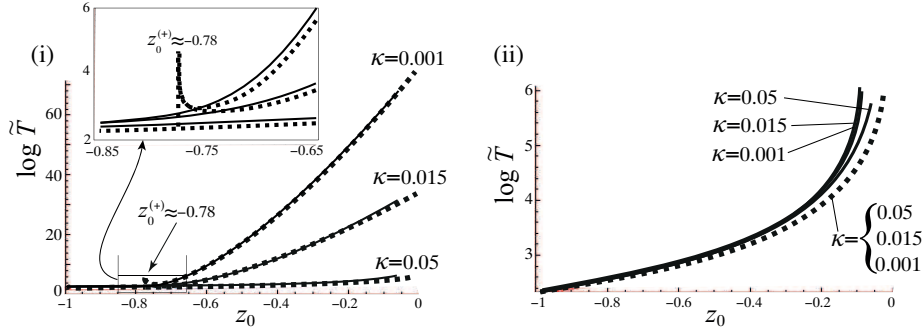}
\caption{\sf\small Plots of $\log\tilde T$, where $T\approx\ee\tilde T$ is the mean escape time to a distance $r=0.1$, for the friction oscillator with smoothing parameter $\ee=0.01$ and friction strength $\alpha=1$. The smoothing function is the cubic polynomial in (\ref{eq:phifric0}) with: (i) $\mu=3$, (ii) $\mu=\frac{1}{2}$. The solid curves show numerical evaluations of the exact integral (\ref{eq:Ttilde}), the dotted curves show the asymptotic value (\ref{eq:Ttildesmall}). Plots are made for different noise amplitudes $\kappa$ as labelled. In (iii) some of the curves, particularly the asymptotic values, are indistinguishable. Inset in (i): a magnification of the Stokes discontinuity at the fold point $z_0=z_0^{(+)}$; in (ii) there is no Stokes discontinuity. 
\label{fig:fric1}
}
\end{center}
\end{figure}

\begin{figure}[h!]
\begin{center}
\includegraphics[width=0.6\textwidth]{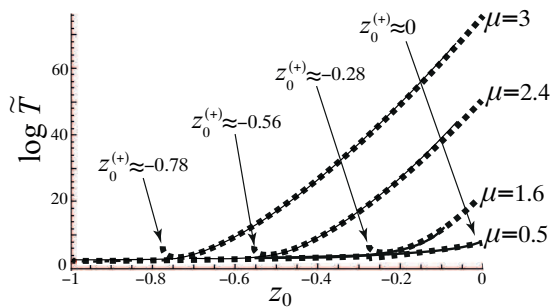}
\caption{\sf\small Plots of the exact integral (solid curves) and asymptotic value (dotted curves) for $\log\tilde T$, as in \fig{fric1}. Each curve is for a fixed noise amplitude $\kappa=0.01$, with $\ee=0.01$ and $\alpha=1$. Four different values of $\mu$ are shown as labelled. The Stokes discontinuities at $z_0=z_0^{(+)}$ are indicated for each case. 
\label{fig:fric2}
}
\end{center}
\end{figure}

By simply substituting (\ref{fr1})-(\ref{fr2}) into (\ref{eq:Ttilde}) and (\ref{eq:Ttildesmall}) we obtain, respectively, the exact and asymptotic values of the mean escape time $T\approx\ee\tilde T$. In \fig{fric1} and \fig{fric2} we plot these for different values of the noise amplitude $\kappa$ and smoothing parameter $\mu$, fixing the size of the smoothing region as $\ee=0.01$, the friction strength as $\alpha=1$, and considering escape to a distance $r=0.1$. 
With these values, the dynamics found on $\sp = 0$ corresponds to:
\begin{itemize}
\item[-] sliding (\fig{flows}(i)) for $0 < z < 2$, 
\item[-] crossing (\fig{flows}(ii.a)) for $z < z_0^{(+)}$, 
\item[-] spurious~sliding (\fig{flows}(ii.b)) for $z_0^{(+)} < z < 0$ .
\end{itemize}

\fig{fric1} shows the effect of different noise amplitudes. In (i), for which $\mu=3$, the mean escape time is strongly dependent on the noise amplitude $\kappa$.
For smaller noise amplitudes, and values of $z_0$ approaching $z=0$ (the boundary of sliding region),
$T\approx\ee\tilde T$ scales with $e^{1/\kappa^2}$. 
This outcome is in stark contrast to that predicted when the Filippov convention is applied to the discontinuous system (\ref{eq:osc2})-(\ref{eq:forcePWC}), for which escape to a distance $r$ from $\sp = 0$
occurs in a time $T=O(r)$ for all $z<0$. 
For larger noise amplitudes the mean escape time is small, simply $\tilde T\approx r/z$ throughout the crossing region, consistent with Filippov dynamics. In (ii), with $\mu=1/2$, the escape time is only weakly dependent on the noise amplitude (in fact the asymptotic value is $\kappa$-independent), and is small (compared to values in (i))
and therefore consistent with Filippov dynamics. 

The key features of \fig{fric1} are understood as follows. 
When $\mu=3$ in \fig{fric1}(i), the friction force $\cal F$ has two prominent turning points (see \fig{fricforce}). These create a potential well in the function $V(v)-V(u)$ at the point $(u,v)=(\tilde \x_1,\tilde \x_2)$, which exists for $z_0$
in the range $z_0^{(+)}\approx-0.78<z_0<0$.
For those $z_0$ values where the potential well exists, the first term in (\ref{eq:Ttildesmall}) dominates, so that $\tilde T$ is strongly $\kappa$-dependent. For values of $z_0$ to the left of $z_0^{(+)}$, the potential well disappears and its contribution to (\ref{eq:Ttildesmall}) is eliminated by the Stokes multiplier, leaving only the smaller second term. The switching of the Stokes multiplier creates a discontinuity, shown magnified in the figure. (Stokes' discontinuities are a well understood artifact of the leading order approximation, which can be smoothed by local approximation about $z_0=z_0^{(+)}$, or more powerfully by uniform approximation \cite{berry-uniform} for arbitrary $z_0$).
In \fig{fric1}(ii) for which $\mu=1/2$,
the friction force has no turning points, so the function $V(v)-V(u)$ has no potential well, therefore the Stokes multiplier is zero everywhere, and the escape time is dominated by the $\kappa$-independent second term in (\ref{eq:Ttildesmall}). 

\fig{fric2} shows four simulations with a small noise amplitude $\kappa=0.1$, for three different values of the parameter $\mu$, in which the potential well in $V(v)-V(u)$ exists for $z_0^{(+)}<z_0<0$, leading to exponentially large escape times to the right of a Stokes discontinuity at $z_0=z_0^{(+)}$, and small escape times to the left. 

In both \fig{fric1} and \fig{fric2}, the exact and asymptotic calculations are shown to be in close agreement for $\ee=0.01$. More generally, of course, the accuracy of the asymptotic approximation depends on the orders of $\ee$ and $r$, as indicated in (\ref{eq:Ttildesmall}). 

More importantly, 
by (\ref{eq:Ttildesmall}) we conclude that
exponentially large escape times occur only when the noise amplitude $\kappa$ and smoothing stiffness $\ee$ satisfy
\begin{equation}
\kappa\gtrsim\sqrt{-\ee/\log(\ee)}\approx0.05 \;.
\end{equation}
This improves on an estimate of $\kappa\gtrsim\ee$ made in \cite{j13error}, which heuristically considered general sources of error of size $\kappa$, instead of the specific white noise considered here. 


\section{Discussion}
\label{sec:CONC}
\setcounter{equation}{0}

We have shown rigoroualy that smoothing a discontinuous system can lead to dynamics that is starkly at odds with what Filippov's method would predict from the discontinuous system. 
For instance, solutions to a smoothed system may remain trapped near a switching surface
when Filippov's convention indicates that solutions to the original system head away from the surface.
However, if the smoothed system is also subject to noise,
the resulting stochastic dynamics may be practically independent of the nature of the smoothing function.
If the noise amplitude is sufficiently large it dominates the deterministic smoothing.
Indeed, in many applications
additional modelling assumptions may have the effect of smoothing a switching surface
out over a range of phase space that is small relative to uncertainties and randomness.
We can think of the noise as acting as a second level of smoothing.
In this case noise inhibits non-Filippov solutions and the use of Filippov's convention is justified.
This is extremely beneficial because, from a modelling perspective,
one typically has little knowledge about the nature of the smoothing function.

In order to quantify this behaviour we introduced equation (\ref{eq:genFilippov})
to describe a general Filippov system near a smooth switching surface.
The stochastic differential equation (\ref{eq:genSDE})
is the result of smoothing (\ref{eq:genFilippov}) within a distance $\ee$ of the switching surface
and adding white Gaussian noise of amplitude $\kappa$.
We used asymptotic methods to study forward evolution from a point on the switching surface
over a short time-frame.

If the noise amplitude $\kappa$ is large relative to the range of the smoothing $\ee$, specifically $\kappa \gg \sqrt{\ee}$,
then the solutions to (\ref{eq:genSDE}) behave like Filippov solutions to (\ref{eq:genFilippov}).
When the initial point lies on an attracting sliding region,
the fraction of time spent within a distance $\ee$ of the switching surface
is small (\ref{eq:Plarge}),
and therefore the effect of the smoothing function $F$ on the dynamics may be neglected.
When the initial point lies on a crossing region,
to leading order, the mean time taken for solutions
to escape a small neighbourhood of the switching surface
is independent of $F$ (\ref{eq:Tlarge}).
Overall, in this scenario the noise suppresses dynamics generated by $F$.

Alternatively, if $\kappa \ll \sqrt{\ee}$
then the nature of the smoothing function $F$ may have a significant effect on the dynamics.
When the initial point lies on an attracting sliding region,
the fraction of time spent within a distance $\ee$ of the switching surface
is near $1$ (\ref{eq:Psmall}).
Consequently, forward evolution is primarily determined by $F$.
In the case of crossing, if the smoothing creates a potential well near the switching surface,
then dynamics may become stuck in the well with high probability (\ref{eq:Tsmall}).
Overall, in this scenario the behaviour of $F$ is important to the dynamics.

We expect that the effect of smoothing and noise on forward evolution from a repelling sliding region (see \cite{jh09,j11prl})
will be analogous to that for crossing, except that solutions may escape
a neighbourhood of the switching surface on both sides.
Furthermore, the smoothing function $F$ may generate a potential well
regardless of whether we are considering a region of crossing, attracting sliding, or repelling sliding.

In \S\ref{sec:FRICTION} we applied these results to the dry-friction oscillator. For a Coulomb-like piecewise-smooth model with a step-function in the friction force, Filippov's convention for sliding and crossing dynamics can be applied. We showed that forming a continuous friction model with a cubic characteristic leads to very different dynamics. Via the stochastic analysis the cubic function was shown to give rise to a potential well, where the system could linger at the discontinuity and thus exhibit spurious sticking in regions where Coulomb's (and Filippov's) model would predict immediate passage through the switching surface. Such spurious dynamics fits qualitatively with the experimental observation that static and kinetic friction coefficients are often unequal \cite{krim,popp,wojewoda08}. 

The particular stochastic formulation (\ref{eq:genSDE}) was, in part, chosen for simplicity.
If the white noise is substituted for coloured noise (which may be more physically realistic)
then the reduction of the system to one dimension that was achieved in \S\ref{sec:SCALE} is not possible.
However, basic Monte-Carlo simulations reveal that coloured noise
induces the same qualitative behaviour when the correlation time is short.
The effects of coloured noise with long correlation times
on dynamics near the switching surface remain to be investigated.
One may also consider what happens if noise is enhanced near the switching surface, perhaps because it stems from the discontinuity itself, as with measurements taken by a relay control system.
It remains to determine whether such ``state noise'' will typically dominate the influence of the smoothing function.


\appendix

\section{Supplementary calculations}

\subsection{Derivation of (\ref{eq:Ttilde})}
\label{sec:TDERIV}
\setcounter{equation}{0}

The probability that a solution to (\ref{eq:reducedSDEtilde}), starting from $\tilde{y} = 0$,
eventually escapes $[-\tilde{r},\tilde{r}]$, is $1$.
Consequently integration of (\ref{eq:FPE}) with respect to $\tilde{t}$ over $[0,\infty)$ yields
\begin{equation}
-\delta(\tilde{\x}) = \frac{d}{d\tilde{\x}}
\left( -\phi R + \frac{\tilde{\kappa}^2}{2} \frac{dR}{d\tilde{\x}} \right) \;,
\label{eq:odeR}
\end{equation}
where we have let $R(\tilde{\x}) = \int_0^\infty p(\tilde{\x},\tilde{t}) \,d\tilde{t}$.
To solve (\ref{eq:odeR}) subject to the boundary conditions $R(\pm \tilde{r}) = 0$,
we first integrate (\ref{eq:odeR}) to obtain
\begin{equation}
-\phi R + \frac{\tilde{\kappa}^2}{2} \frac{dR}{d\tilde{\x}} = C - H(\tilde{\x}) \;,
\end{equation}
where $H$ is the Heaviside function and $C$ is a constant.
Through the use of an integrating factor, further integration produces
\begin{equation}
R(\tilde{\x}) = \frac{2}{\tilde{\kappa}^2} \,{\rm e}^{\frac{-2 V(\tilde{\x})}{\tilde{\kappa}^2}}
\int_{\tilde{\x}}^{\tilde{r}} \left( H(v) - C \right) {\rm e}^{\frac{2 V(v)}{\tilde{\kappa}^2}} \,dv \;,
\end{equation}
where we have chosen the limits of integration to automatically satisfy $R(\tilde{r}) = 0$.
The requirement $R(-\tilde{r}) = 0$ implies that $C$ is given by (\ref{eq:C}).
By (\ref{eq:Ttilde0}) we have $\tilde{\T} = \int_{-\tilde{r}}^{\tilde{r}} R(u) \,du$,
and the expression (\ref{eq:Ttilde}) follows from reversing the order of integration. 

\subsection{Bounding $C$ for (\ref{eq:TtildeCzero})}\label{sec:C}

Here we derive an upper bound for $C$ (\ref{eq:C}).
Given any $\Delta > 0$, we have
\begin{equation}
\int_{-\tilde{r}}^0 {\rm e}^{\frac{2 V(w)}{\tilde{\kappa}^2}} \,dw
\ge \int_{-\tilde{r}}^{-1} {\rm e}^{\frac{2 V(w)}{\tilde{\kappa}^2}} \,dw
= \frac{\tilde{\kappa}^2}{2 a^-} \left( {\rm e}^{\frac{2 a^- (\tilde{r} - 1)}{\tilde{\kappa}^2}} - 1 \right)
\ge \frac{\tilde{\kappa}^2 (1-\Delta)}{2 a^-} \,{\rm e}^{\frac{2 a^- \tilde{r} (1-\Delta)}{\tilde{\kappa}^2}} \;,
\end{equation}
for sufficiently small $\ee$, $\kappa$ and $r$.
We let $V_{\rm max} = {\rm max}_{w \in [0,1]} V(w)$.
Then
\begin{equation}
\int_0^{\tilde{r}} {\rm e}^{\frac{2 V(w)}{\tilde{\kappa}^2}} \,dw
\le {\rm e}^{\frac{2 V_{\rm max}}{\tilde{\kappa}^2}}
\left( 1 + \int_1^{\tilde{r}} {\rm e}^{\frac{-2 a^+ (w-1)}{\tilde{\kappa}^2}} \,dw \right)
\le \left( 1 + \frac{\tilde{\kappa}^2}{2 a^+} \right) {\rm e}^{\frac{2 V_{\rm max}}{\tilde{\kappa}^2}} \;.
\end{equation}
From (\ref{eq:C}), we have
\begin{align}
C &= \frac{1}{1 + \frac{\int_{-\tilde{r}}^0 {\rm e}^{\frac{2 V(w)}{\tilde{\kappa}^2}} \,dw}
{\int_0^{\tilde{r}} {\rm e}^{\frac{2 V(w)}{\tilde{\kappa}^2}} \,dw}}
\le \frac{1}{1 + \frac{\frac{\tilde{\kappa}^2 (1-\Delta)}{2 a^-} \,{\rm e}^{\frac{2 a^- \tilde{r} (1-\Delta)}{\tilde{\kappa}^2}}}
{\left( 1 + \frac{\tilde{\kappa}^2}{2 a^+} \right) {\rm e}^{\frac{2 V_{\rm max}}{\tilde{\kappa}^2}}}} \nonumber \\
&\le \frac{(1+\Delta) a^-}{1-\Delta} \left( \frac{1}{a^+} + \frac{2}{\tilde{\kappa}^2} \right)
{\rm e}^{\frac{-2 a^- \tilde{r} \left( 1-\Delta-\frac{V_{\rm max}}{\tilde{r}} \right)}{\tilde{\kappa}^2}} \nonumber \\
&\le \left( 1 + 4 \Delta \right) a^- \left( \frac{1}{a^+} + \frac{2}{\tilde{\kappa}^2} \right)
{\rm e}^{\frac{-2 a^- \tilde{r} \left( 1-2\Delta \right)}{\tilde{\kappa}^2}} \;,
\end{align}
where, in the last inequality we have assumed $\tilde{r}$ is sufficiently large that
$\frac{V_{\rm max}}{\tilde{r}} \le \Delta$.
Substituting $\Delta = \frac{1}{4}$ gives the result in the text.

\subsection{The double integral (\ref{eq:secondTerm})}\label{sec:doubleint}

The double integral is independent of $\frac{1}{\tilde{r}}$,
so we simply evaluate it asymptotically in $\tilde{\kappa}$.
Since the maximum of the exponent is attained at $(u,v) = (\tilde{\x}_1,\tilde{\x}_2)$,
we employ the integral substitution:
\begin{equation}
\hat{u} = \frac{u - \tilde{\x}_1}{\tilde{\kappa}} \qquad
\hat{v} = \frac{v - \tilde{\x}_2}{\tilde{\kappa}} \;.
\end{equation}
This yields
\begin{align}
&\frac{2}{\tilde{\kappa}^2} \int_0^1 {\rm e}^{\frac{2 V(v)}{\tilde{\kappa}^2}}
\int_{-1}^v {\rm e}^{\frac{-2 V(u)}{\tilde{\kappa}^2}} \,du \,dv \nonumber \\
&= 2 \int_{\frac{-\tilde{\x}_2}{\tilde{\kappa}}}^{\frac{1 - \tilde{\x}_2}{\tilde{\kappa}}}
\int_{\frac{-1 - \tilde{\x}_1}{\tilde{\kappa}}}^{\hat{v} + \frac{\tilde{\x}_2 - \tilde{\x}_1}{\tilde{\kappa}}}
{\rm e}^{\frac{-2}{\tilde{\kappa}^2} \left( V(\tilde{\x}_1) +
\frac{1}{2} V''(\tilde{\x}_1) \tilde{\kappa}^2 \hat{u}^2 - V(\tilde{\x}_2) -
\frac{1}{2} V''(\tilde{\x}_2) \tilde{\kappa}^2 \hat{v}^2 + O(\tilde{\kappa}^3) \right)}
\,d\hat{u} \,d\hat{v} \;.
\end{align}
Laplace's method justifies taking the limits of integration to $\pm \infty$,
with which explicit integration produces
\begin{equation}
\frac{2}{\tilde{\kappa}^2} \int_0^1 {\rm e}^{\frac{2 V(v)}{\tilde{\kappa}^2}}
\int_{-1}^v {\rm e}^{\frac{-2 V(u)}{\tilde{\kappa}^2}} \,du \,dv
= 2 {\rm e}^{\frac{2}{\tilde{\kappa}^2} \left( V(\tilde{\x}_2) - V(\tilde{\x}_1) \right)}
\left( \frac{\sqrt{\pi}}{\sqrt{V''(\tilde{\x}_1)}}
\frac{\sqrt{\pi}}{\sqrt{-V''(\tilde{\x}_2)}} + O(\tilde{\kappa}) \right) \;.
\end{equation}
Note $V''(\tilde{\x}) = -A'(\tilde{\x})$.

\end{document}